\def\draft{n}
\documentclass{amsart}
\usepackage{fullpage,amssymb,epic,eepic,epsfig,amscd,pb-diagram,lamsarrow,pb-lams}

\theoremstyle{plain}

\newtheorem{theorem}{Theorem}
\newtheorem{proposition}{Proposition}[section]
\newtheorem{lemma}[proposition]{Lemma}
\newtheorem{corollary}[proposition]{Corollary}

\theoremstyle{definition}
\newtheorem{definition}[proposition]{Definition}

\theoremstyle{remark}

\newtheorem{remark}[proposition]{Remark}

\def\printname#1{
	\if\draft y
		\smash{\makebox[0pt]{\hspace{-0.5in}
			\raisebox{8pt}{\tt\tiny #1}}}
	\fi
}

\newcommand{\psdraw}[2]
         {\begin{array}{c} \hspace{-1.3mm}
	\raisebox{-4pt}{\epsfig{figure=draws/#1.eps,width=#2}}
	\hspace{-1.9mm}\end{array}}

\newlength{\standardunitlength}
\catcode`\@=11
\long\def\@makecaption#1#2{%
     \vskip 10pt

\setbox\@tempboxa\hbox{
       \small\sf{\bfcaptionfont #1. }\ignorespaces #2}%
     \ifdim \wd\@tempboxa >\captionwidth {%
         \rightskip=\@captionmargin\leftskip=\@captionmargin
         \unhbox\@tempboxa\par}%
       \else
         \hbox to\hsize{\hfil\box\@tempboxa\hfil}%
     \fi}
\font\bfcaptionfont=cmssbx10 scaled \magstephalf
\newdimen\@captionmargin\@captionmargin=2\parindent
\newdimen\captionwidth\captionwidth=\hsize
\catcode`\@=12

\newcommand{\tr}{\operatorname{tr}}

\def\lbl#1{\label{#1}\printname{#1}}

\def\BZ{\mathbb Z}
\def\BQ{\mathbb Q}
\def\BR{\mathbb R}
\def\BC{\mathbb C}

\def\A{\mathcal A}
\def\B{\mathcal B}

\def\D{\Delta}

\def\a{\alpha}

\def\La{\Lambda}
\def\l{\lambda}

\def\ihs{integral homology 3-sphere}

\def\la{\langle}
\def\ra{\rangle}

\def\sub{\subseteq}

\def\strutb#1#2#3{\overset{#1}{\underset{#2}{ 
\begin{array}{c} \vspace{0.0cm}
\uparrow 
\vspace{-0.25cm} \\        
| \vspace{-0.45cm} \\      
\bullet \vspace{0.00cm}   
\end{array} }}\! #3}

\def\AS{\mathrm{AS}}
\def\IHX{\mathrm{IHX}}

\def\Lhat{\hat\Lambda}
\def\Lloc{\Lambda_{\mathrm{loc}}}

\def\loc{\mathrm{loc}}

\def\Zrat{Z^{\mathrm{rat}}}

\def\Aut{\mathrm{Aut}}
\def\hair{\mathrm{Hair}}
\def\hairnu{\mathrm{Hair}^{\nu}}

\def\alg{\mathrm{alg}}
\def\fg{\mathfrak{g}}

\def\ft{\mathfrak{t}}

\def\Ad{\mathrm{Ad}}
\def\ad{\mathrm{ad}}

\def\longto{\longrightarrow}

\begin{document}


\title[Beads: From Lie algebras to Lie groups]{Beads: From Lie algebras to 
Lie groups}

\address{School of Mathematics \\
         Georgia Institute of Technology \\
         Atlanta, GA 30332-0160, USA}
\email{stavros@math.gatech.edu}

\thanks{Partially supported by NSF and BSF. \\
        This and related preprints can also be obtained at
{\tt http://www.math.gatech.edu/$\sim$stavros } 
\newline
1991 {\em Mathematics Classification.} Primary 57N10. Secondary 57M25.
\newline
{\em Key words and phrases: Kontsevich integral, rational form, weight 
systems, Lie algebras, Lie groups, beads.} 
}

\date{
This edition: March 1, 2003 \hspace{0.5cm} First edition: January 2, 2002.}


\begin{abstract}
The Kontsevich integral of a knot is a powerful invariant which takes values
in an algebra of trivalent graphs with legs. Given a Lie algebra, the
Kontsevich integral determines an invariant of knots (the so-called colored 
Jones function) with values in the symmetric algebra of the Lie
algebra. Recently A. Kricker and the author constructed a rational form of
the Kontsevich integral which takes values in an algebra of trivalent
graphs with beads. After replacing beads by an exponential legs, this 
rational form recovers
the Kontsevich integral. The goal of the paper is to explain the relation
between beads and functions defined on a Lie group. As an application,
we provide a rational form for the colored Jones function of a knot,
conjectured by Rozansky.
\end{abstract}

\maketitle



\section{Introduction}
\lbl{sec.intro}
\subsection{Finite type invariants of knots, Feynmann diagrams and Lie 
algebras}
\lbl{sub.history}

A {\em knot} is a (smooth) imbedding of a circle $S^1$ in 3-space. An 
{\em immersed knot} is a (smooth) immersion of $S^1$ in 3-space with
transverse double points. Any numerical invariant $f$ of knots which
takes values in an abelian group $A$ can be extended to an invariant
of immersed knots via the rule:
$$
f\left(\psdraw{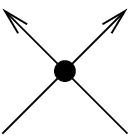}{.3in}\right)=
f\left(\psdraw{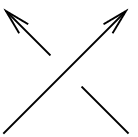}{.3in}\right)-f\left(\psdraw{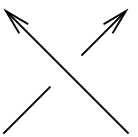}{.3in}\right)
$$
In the early nineties, Vassiliev (and independently, Goussarov) introduced
the notion of a {\em finite type invariant} of knots, that is of a function
$f: \mathrm{Knots}\longto A$, which, when extended to a function on the
set of immersed knots, vanishes on immersed knots with sufficiently many
double points.
Shortly afterwards, Kontsevich constructed an invariant 
$$
Z: \mathrm{Knots} \longto \A(\star)
$$
(the so-called {\em Kontsevich integral}) where 
$\A(\star)$ is the vector space over the rational numbers $\BQ$ generated
by unitrivalent graphs (with vertex orientations), modulo the well-known
antisymmetry $\AS$ and $\IHX$ relations. The Kontsevich integral $Z$
has two key properties:
\begin{itemize}
\item
$Z$ is a {\em universal} finite type invariant. That is, it determines
every $\BQ$-valued finite type invariant of knots.
\item
$Z$ determines the {\em Jones polynomial} of a knot, and more generally
the invariants of knot that are defined using {\em quantum groups} and their
representations.
\end{itemize}
Let us briefly recall the latter property. 
Given a Lie algebra $\fg$ with an invariant inner product, there 
is a map (often called a {\em weight system})
\begin{equation}
\lbl{eq.weight}
W_{\fg}^h : \A(\star) \longto S(\fg)^{\fg}[[h]], 
\end{equation}
where $S(\fg)$ is the {\em symmetric} algebra of the vector space $\fg$.
\begin{definition}
\lbl{def.colored}
We will call the map 
$$
W_{\fg}^h \circ Z: \mathrm{Knots} \longto S(\fg)^{\fg}[[h]]
$$
the $\fg$-{\em colored Jones function} of a knot, and we will denote
it by $J_{\fg}$.
\end{definition}
In a sense, the colored Jones function of a knot is a generating function
of the quantum group invariants of a knot. More precisely, given 
an irreducible representation $V$ of $\fg$, its 
evaluation at its dominant weight gives rise to a linear map
$$
S(\fg)^{\fg}[[h]] \longto \BQ[[h]]
$$
The image of the Kontsevich integral under the composition $W_{\fg,V}$ of the
above two maps is an element of the ring $\BZ[q^{\pm 1}]$ (where $q=e^h$)
and coincides with the quantum group invariant of knots, using the $(\fg,V)$
data. If $\fg=\mathfrak{sl}_2$ and $V=\BC^2$ is the defining representation
of $\mathfrak{sl}_2$, then $W_{\mathfrak{sl}_2, \BC^2} \circ Z$ coincides
with the Jones polynomial.

\subsection{A rational form of the Kontsevich integral, Feynmann diagrams 
and beads}
\lbl{sub.rational}

Recently, A. Kricker and the author \cite{GK} constructed a {\em
rational form} $\Zrat$ of the Kontsevich integral of a knot
$$
\Zrat: \mathrm{Knots} \longto 
\A^0(\Lloc):=\B(\La\to\BZ) \times \A(\Lloc)
$$
which consists of a 'matrix part' and a 'graph-part', where 
\begin{itemize}
\item
$\La=\BZ[t,t^{-1}]$, $\Lloc$ is the localization of $\La$ with respect
to elements $f \in \La$ such that $f(1)=1$.
\item
$\B(\La\to\BZ)$ 
is a quotient of the set of Hermitian matrices over 
$\La=\BZ[t,t^{-1}]$ which are invertible over $\BZ$, modulo the equivalence
$A \sim B$ iff $A \oplus D=P(B\oplus E)P^\star$ for diagonal matrices $D,E$
with monomials in $t$ on the diagonal and for $P$ invertible over $\La$.
\item
$\A(\Lloc)$ is the completed vector space of trivalent graphs with oriented 
vertices and edges, and with elements (often called {\em beads}) of $\Lloc$ 
associated to each edge, modulo the $\AS$, $\IHX$, Orientation Reversal,
Linearity and Holonomy relations shown in Figure \ref{relations4}. 
$\A(\Lloc)$ is a graded vector space, where the degree of a diagram is
half the number of trivalent vertices, and the completion refers to the above
grading.
\end{itemize}

\begin{figure}[htpb]
$$ 
\psdraw{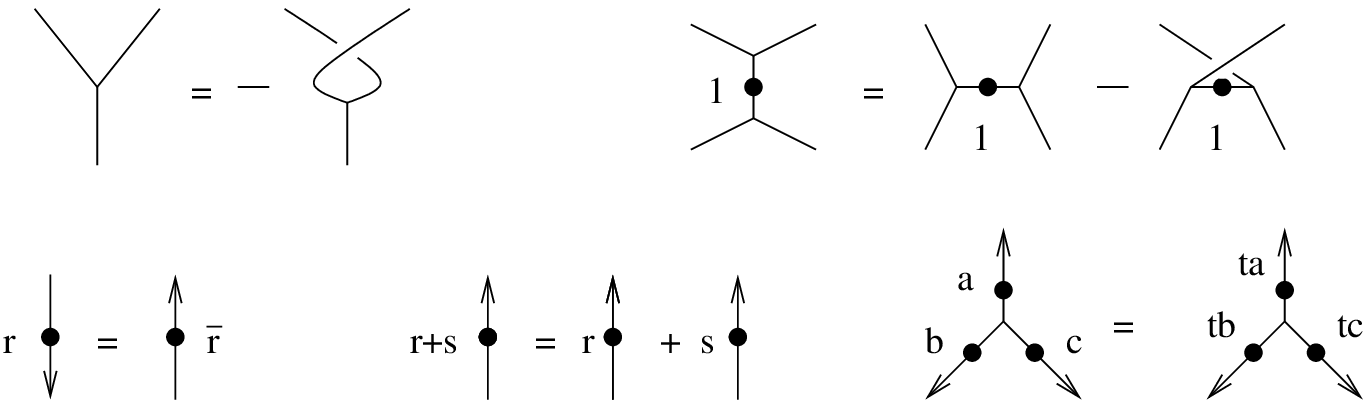}{4.5in} 
$$
\caption{The $\AS$, $\IHX$ (for arbitrary orientations of the edges),
Orientation Reversal, Linearity and Holonomy Relations.}\lbl{relations4}
\end{figure}

The rational form $\Zrat$ has two key properties:
\begin{itemize}
\item
It is a {\em universal finite type invariant of knots}, where the analogue
of crossing change (in the definition of finite type invariants) is
replaced by the {\em null move} described in terms of 
surgery on claspers with nullhomologous leaves in a knot complement.
For details, see \cite{GR} and \cite{GK}.
\item
$\Zrat$ determines the Kontsevich integral of a knot, after we {\em replace
beads by hair}.
\end{itemize}
Let us explain the last phrase. Given a trivalent graph $s$ with beads in 
$\La$, let
$\hair(s) \in \A(\star)$ denote the result of replacing each bead $t$
by an exponential of hair (i.e., legs):
$$
\strutb{}{}{t} \to \sum_{n=0}^\infty \frac{1}{n!} \printname{attacchn}
	\setlength{\unitlength}{0.03\standardunitlength}
	\begin{array}{c}  \hspace{-1.7mm}
         	\raisebox{-8pt}{
\begingroup\makeatletter\ifx\SetFigFont\undefined%
\gdef\SetFigFont#1#2#3#4#5{%
  \reset@font\fontsize{#1}{#2pt}%
  \fontfamily{#3}\fontseries{#4}\fontshape{#5}%
  \selectfont}%
\fi\endgroup%
{\renewcommand{\dashlinestretch}{30}
\begin{picture}(1165,1239)(0,-10)
\path(12,12)(12,1212)
\path(42.000,1092.000)(12.000,1212.000)(-18.000,1092.000)
\path(12,912)(312,912)
\path(192.000,882.000)(312.000,912.000)(192.000,942.000)
\path(12,612)(312,612)
\path(192.000,582.000)(312.000,612.000)(192.000,642.000)
\path(12,312)(312,312)
\path(192.000,282.000)(312.000,312.000)(192.000,342.000)
\put(612,537){\makebox(0,0)[lb]{$n$ legs}}
\end{picture}
} }
         	\hspace{-1.9mm}
	\end{array}
 
$$
This map can be defined even if the beads are elements of $\Lhat=\BQ[[t-1]]$,
and in particular if the beads are elements of $\Lloc \sub \Lhat$. Notice
that $\hair(s)$ is a series of unitrivalent graphs that contain no wheels,
where, for example, a {\em wheel} with $4$ hair is:
$$
\psdraw{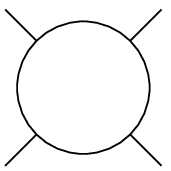}{0.2in}
$$
Next, we add wheels from the matrix part of $\A^0(\Lloc)$ as follows. 
We define 
$$
\hairnu: \A^0(\Lloc)\longto \A(\star)
$$
by 
$$
(A,s) \in \B(\La\to\BZ) \times \A(\Lloc) \mapsto
\nu \, \exp\left(-\frac{1}{2} \tr\log(A)(e^h)|_{h^n \to
w_n} \right) \, \hair(s)
$$
where $\nu=Z(\text{unknot})$, $w_n$ is the wheel with $n$ legs 
and we think of $\A(\star)$ as an algebra with the disjoint union
multiplication of graphs.

\cite[Theorem 1.3]{GK} states that
\begin{equation}
\lbl{eq.ZratZ}
Z=\hairnu \circ \Zrat .
\end{equation}

Let us end this introduction to the $\Zrat$ invariant with three comments. 
For knots with a fixed Alexander polynomial $\D$, the graph-part of the
$\Zrat$ invariant takes values in graphs with beads in the abelian subgroup
$$
\La_{\D}=\frac{1}{\D} \La
$$
of $\La$. Furthermore, if we restrict $\Zrat$ to
knots with trivial Alexander polynomial, then the matrix part of $\Zrat$
is trivial and there is no need to consider graphs with beads in $\Lloc$.
Namely, we have:
$$
\Zrat: \text{Knots with trivial Alexander polynomial} \longto
\A(\La).
$$
As a first approximation to understanding of the rational form $\Zrat$,
the reader may restrict their attention to knots with trivial Alexander
polynomial.

A last comment: the Kontsevich integral $Z$ and its rational form $\Zrat$
can be defined for knots $K$ in \ihs s $M$, i.e., manifolds $M$ such that
$H_\star(M,\BZ)=H_\star(S^3,\BZ)$.

\subsection{Weight systems for diagrams with beads}
\lbl{sub.wbeads}

The paper is concerned with defining a concept of a weight system
for diagrams with beads. As an application, we will deduce a rational
form for the colored Jones function of a knot. As we will see, the
analogue of the weight system map \eqref{eq.weight} uses the Lie group,
rather than its Lie algebra. This explains the title of the paper.

For a reference on Lie groups and their representations, see \cite{CSM}.
Consider a compact connected (not necessarily simply connected)
Lie group $G$ with semisimple Lie algebra 
$\fg_{\BR}$, complexification $\fg=\fg_{\BR}\otimes \BC$, 
and let $C(G)$ the algebra of continuous complex-valued functions on $G$. 
$G$ acts on $C(G)$ by conjugation. Consider the
subspace $C_{\alg}(G)$ of {\em almost invariant functions}, that is functions 
$f \in C(G)$ whose image under $G$ lies in a finite dimensional subspace of 
$C(G)$. It turns out that $C_{\alg}(G)$ is a subalgebra of $C(G)$. Since 
$C_{\alg}(G)$ is perhaps not a familiar enough object to a topologist, we 
discuss several alternative views of it in the following:

\begin{lemma}
\lbl{lem.views}
There are canonical isomorphisms between $C_{alg}(G)$ and 
\begin{itemize}
\item[(a)]
$\oplus_{E \in \hat{G}} E' \otimes E$, where $\hat G$ is the set of
unitary irreducible representations of $G$.
\item[(b)]
The ring of {\em representative functions} on $G$, 
where a representative
function on $G$ is a function $G \to \BC$ of the form
$$
f_{M,\xi,\eta}(g)=\la \eta, g \xi \ra
$$
where $M$ is a finite dimensional unitary representation of $G$ and $\xi,\eta
\in M$. In other words, $f_{M,\xi,\eta}(g)$ is a matrix element of the action
of $g$ on $M$ with respect to a suitable basis.
\item[(c)]
The coordinate ring of the linear algebraic group $G_{\BC}$.
\end{itemize}
\end{lemma}

\begin{proof}
All these are restatements of the {\em Peter-Weyl} theorem. For an excellent
discussion of these ideas, we refer the reader to \cite{CSM}. For (a) see 
\cite[pp.91-100]{CSM}. For (b) and (c), see \cite[p.92]{CSM}.
\end{proof}

Let $R(G)$ denote the algebra of {\em characters} of unitary representations 
of $G$, a subalgebra of the {\em class functions} $C(G)^G$, and let
$R(T)^W$ denote the algebra of {\em Weyl-invariant} characters of a 
{\em maximal torus} $T$ of $G$.

\begin{corollary}
\lbl{cor.1}
We have isomorphisms of algebras:
$$
C_{\alg}(G)^G \cong R(G) \cong R(T)^W.
$$
\end{corollary}

\begin{proof}
Characters are invariant functions on $G$, thus there is an inclusion
$R(G) \to C_{\alg}(G)^G$. The first description of $C_{\alg}(G)$ above and 
Shur's lemma (in the form $(E' \otimes E)^G\cong \BC$) implies that
$C_{\alg}(G)^G \cong \oplus_{E \in \hat{G}} \BC$, from which follows that
the inclusion $R(G) \to C_{\alg}(G)^G$ is an isomorphism.

The restriction of characters from $G$ to $T$ induces an isomorphism 
$R(G)\cong R(T)^W$. 
This is the content of Chevalley's theorem \cite[Section 23.1]{Hu} and also
\cite[Theorem 7.28]{BGV}, and is a consequence of the fact that given 
$g \in G$, there exist $h \in G$ and $t \in T$ unique up to the action of 
$W=N_G(T)/T$ (where $N_G(T)$ is the normalizer of $T$ in $G$) such that 
$g=h^{-1}th$. 
\end{proof}

We will need a subspace of the algebra $C_{\alg}(G)$ of almost invariant
functions. Given a finite set $S$ of conjugacy classes in $G$, let
$R_S(G)$ denote the class functions on $G-S$, and $R_R(T)^W$ denote the
$W$-invariant characters on $T-S$. We define $C_{\alg, S}(G)^G=R_S(G)$. 

The above corollary implies that
\begin{corollary}
\lbl{cor.1S}
For every $S$ as above, we have vector space isomorphisms:
$$
C_{\alg,S}(G)^G \cong R_S(G) \cong R_S(T)^W.
$$
\end{corollary}

Given an element $f \in \La$, consider the subset $S$ of elements in $G$ such 
that $\det(f(\Ad(g)))=0$. It is easy to see that $S$ is a finite union
of conjugacy classes, and that $S$ is empty if $f$ has {\em no roots} on
the unit circle $S^1$.  We will denote the corresponding vector spaces by
$C_{\alg,f}(G)^G$, $R_f(G)$ and $R_f(T)^W$ respectively.

\begin{remark}
\lbl{rem.generators}
For computational reasons, it is useful to know the structure of the algebra 
$R(T)^W$. $R(T)$ can be identified with the group-ring $\BC[\La_w]$ of
the {\em weight lattice} $\La_w \subset \ft^{\star}$
$$
R(T) \cong \BC[\La_w].
$$
For $\l \in \La_w$, we will denote by $e_{\l}$ the associated element of
$\BC[\La_w]$ and $R(T)$. Notice that $e_{\l} e_{\l'}=e_{\l + \l'}$.
The corresponding function
$$
e_{\l}:T \to \BC
$$ 
is given by $e_{\l}(e^x)=e^{\l(x)}$ for $x \in \ft$, the Lie algebra of $T$.
Although we will not use, a theorem of Chevalley states that $R(T)^W$ is 
freely generated by polynomials in the variables $e_{\l}$. 
\qed
\end{remark}

We now have the preliminaries to define a notion of weight system for
diagrams with beads. 

\begin{definition}
\lbl{def.weightb}
Given $G$ as above, define a map
\begin{equation}
\lbl{eq.weightb}
W_G^h:\A(\La_{\D}) \longto C_{\alg,\D}(G)^G[[h]]
\end{equation}
as follows: fix an element $g \in G$, and a trivalent graph $\Gamma$
with beads. Cut $\Gamma$ into a union of $Y$ graphs together with oriented 
edge segments which contain the beads. Color each of the $Y$ graphs with 
elements of a basis $e_a$ of the Lie algebra $\fg$ of $G$, and replace an 
oriented segment with a bead $f(t) \in 1/\D \La$ with $f(\Ad(g)) \in 
\Aut(\fg)$. This can be done as long as $\det(\D(\Ad(g))) \neq 0$, i.e.,
as long as $g \in G-\D$, using the above notation. Then, contract 
the indices, using structure constants on the $Y$ graphs. Varying $g$
defines an element $W_G(\Gamma) \in C(G)$ which is independent of
the basis of $\fg$, and further lies in the subalgebra of $C(G)$ generated
by the entries of matrices of the Adjoint representation of $G$. This implies 
that $W_G(\Gamma) \in C_{\alg}(G)$. It is easy to see that the $\AS, \IHX$ 
and Holonomy relations of $\A(\La)$ are satisfied and that $W_G(\Gamma)$ lies 
in the $G$-invariant part $C_{\alg}(G)^G$ of $C_{\alg}(G)$. Finally, let
$W_G^h=W_G h^{\mathrm{deg}}$, where $\mathrm{deg}$ is half the number of 
trivalent vertices of a diagram.
Next, we extend the above map, namely we define
$$
W_G^h: \A^0(\La_{\D})\longto C_{\alg,\D}(G)^G[[h]]
$$
by 
$$
(A,s) \in \B(\La\to\BZ) \times \A(\La_{\D}) \mapsto
 \, 
\frac{1}{\det( \det(A)|_{t \to Ad(g)})^{1/2}} W_G^h(s) 
$$
\end{definition}

Next we discuss how to compare the weight systems $W_{\fg}^h$ and
$W_G^h$ of Equations \eqref{eq.weight} and \eqref{eq.weightb}. To achieve
this, we need to
discuss the map $\hair$ on the level of Lie algebras. Let us denote by
$$
\exp^{\star}: C_{\alg}(G)^G[[h]]\longrightarrow S(\fg)^{\fg}[[h]]
$$
the map defined by 
$$
\exp^{\star}(f)(\l)=
f(e^{h \l}) \in P(\fg)[[h]] \cong P(\fg^{\star})[[h]]
\cong S(\fg)[[h]]
$$
for $f \in C_{\alg}(G)^G$ and $\l \in \fg$, 
where $P(V)$ denotes the polynomial functions on a vector space $V$
and where $e: \fg \to G$ is the exponential function and where $\fg \cong 
\fg^{\star}$ via an invariant inner product.
Notice the $h$ in the $e^{h \l}$ term, which can be interpreted by
considering the $1$-parameter family $\fg_h$ of Lie algebras on the vector 
space $\fg$ with Lie bracket $[a,b]_h=[a,b]$, where $[\cdot,\cdot]$ is the 
Lie bracket on $\fg$.

Since $\exp^{\star}(f) \in S(\fg)^{\fg}[[h]]$ is invertible when $f(1)=1$, it 
follows that $\exp^{\star}$ descends to a map
$$
\exp^{\star}: C_{\alg,\D}(G)^G[[h]]\longrightarrow S(\fg)^{\fg}[[h]]
$$

\begin{definition}
\lbl{def.Phi}
A slightly normalized version of $\exp^{\star}$ is the map 
$$
\Phi: C_{\alg,\D}(G)^G[[h]]\longrightarrow S(\fg)^{\fg}[[h]]
$$
given by
$$
\Phi(f)(\l)= \frac{j^{1/2}(h \l)}{j^{1/2}(h \rho)} \exp^{\star}(f)(\l).
$$ 
where $j^{1/2}: \fg \longto \BC$ is defined by
$$ j^{1/2}(x) =
  \det\left(\frac{\sinh\ad_{x/2}}{\ad_{x/2}}\right)^{1/2}.
$$
(The square root is well-defined, see \cite{Du}).
\end{definition}
Observe that our map $\Phi$ differs
from the map $\Phi$ of Kashiwara-Verge \cite{KV} by the factor 
$j^{1/2}(h \rho)$, which is independent of a knot. $\Phi$ turns out to 
coincide with the map $\hairnu$
on the level of Lie algebras, as is revealed by the next theorem.

\begin{definition}
\lbl{def.coloredb}
We will call the map 
$$
W_G^h \circ \Zrat \longto C_{\alg,\D}(G)^G[[h]]
$$ 
the $G$-{\em colored Jones function} of a knot, and we will denote it by $J_G$.
\end{definition}
The next theorem, which compares the $G$ and $\fg$-colored Jones function
of a knot, is often called the {\em Rationality Conjecture} for the
$\fg$-colored Jones function. The conjecture was formulated by Rozansky in
\cite{R2} and proven by entirely different means for $\mathfrak{sl}_2$
in \cite{R1}. 

We are now ready to state the main result of the paper.

\begin{theorem}
\lbl{thm.1}
\rm{(a)}
The weight systems $W_{\fg}^h$ and $W_G^h$ 
fit in a commutative diagram:
$$
\begin{diagram}
\node{\A^0(\Lloc)}              
\arrow{e,t}{W_G^h}\arrow{s,r}{\hairnu}
\node{C_{\alg,\loc}(G)^G[[h]]} 
\arrow{s,r}{\Phi} \\
\node{\A(\star)}            
\arrow{e,t}{W_{\fg}^h}
\node{S(\fg)^{\fg}[[h]]}
\end{diagram}
$$
In other words, the image of the map $\hairnu$ coincides with the 
Kashiwara-Vergne map $\Phi$ on the level of Lie algebras. 
\newline
\rm{(b)}
We have: 
$$
J_{\fg}= \Phi \circ J_G.
$$
In other words, if we write
$$
J_G=\sum_{n=0}^\infty Q_{G,n} h^n
$$
for $Q_{G,n} \in C_{\alg,\D}(G)^G[[h]]$, then we have
$$ 
J_{\fg}(\l)= [ \l-\rho] 
\sum_{n=0}^\infty Q_{G,n}|_{e_\mu \to e^{h (\l,\mu)}} h^n \in S(\fg)^{\fg}[[h]]
$$
where $[ \l-\rho]$ is the {\em quantum dimension} of the representation
of $\fg$ of highest weight $\l-\rho$ given by
$$
[\l-\rho] =
\prod_{\a \succ 0} \frac{\sinh h(\a,\l)}{\sinh h (\a,\rho)}.
$$
\end{theorem}

In other words, $J_G$ is equivalent to a sequence $\{Q_{G,n}\}_n$ 
of (partially defined) 
$G$-invariant functions on $G$ associated to a knot. If a knot has trivial
Alexander polynomial, then it follows from Remark \ref{rem.generators}
that $J_G$ is equivalent to a sequence of polynomials, i.e., elements
of $\BC[\La_w]$. It is well known that these polynomials actually lie
in $\BQ[\La_w]$, see also \cite{A}.

At any rate, any evaluation of these sequences of function gives a knot 
invariant. There are two
natural evaluations to consider. Evaluation at $1$ and average over $G$. 
Let us end with a definition. 

\begin{definition}
\lbl{def.average}
Given a Lie group $G$ as above, let us define two invariants
$$
\kappa^h_G: \mathrm{Knots} \longto \BQ[[h]]
$$
$$
\tau^h_G: \text{Knots with trivial Alexander polynomial} \longto \BC[[h]]
$$
by
$$ \kappa^h_G= J_G(1) \hspace{1cm} \text{and} \hspace{1cm}
\tau^h_G= \int_G J_G d\mu
$$
where $\mu$ is the {\em Haar measure} on $G$.
We will call $\tau^h_G$ the {\em Kashaev invariant} of the knot.
\end{definition}

Unpublished work of Thang Le, \cite{Le} 
combined with the results of this paper shows 
that when $G=SU(2)$, the Kasahev invariant is related to the usual Kashaev 
invariant of the knot. More precisely, Le constructs a function
$$
\kappa: \mathrm{Knots} \longto \widehat{\BZ[q]}
$$
where $\widehat{\BZ[q]}=\lim_{\leftarrow n} 
\BZ[q]/((1-q)(1-q^2) \dots (1-q^n))$ is the 
{\em cyclotomic completion} of the ring $\BZ[q]$ such that 
\begin{itemize}
\item
$\kappa$ evaluates to the Kashaev invariants of the knot, i.e., for all
$n$ we have
$$
\kappa(e^{2 \pi i /n})=\bar J_{\mathfrak{sl}_2,\BC^n}(e^{2 \pi i /n})
$$
where $\bar J_{\fg}=J_{\fg}/J_{\fg}(\mathrm{unknot})$.
\item
When composed with the map $\widehat{\BZ[q]} \longto \BQ[[h]]$,
$\kappa$ coincides with our invariant $\kappa^h_{SU(2)}$. 
\end{itemize}

\subsection{Acknowledgement}
The author wishes to thank A. Kricker for numerous conversations and TTQ. Le
for pointing an error in an ealier version of the paper.

\section{Proof of Theorem \ref{thm.1}}
\lbl{sec.weight}

\rm{(a)} It suffices to show that the following diagrams commute:
$$
\begin{diagram}
\node{\A(\La_{\D})}              
\arrow{e,t}{W_G^h}\arrow{s,r}{\hairnu}
\node{C_{\alg,\D}(G)^G[[h]]} 
\arrow{s,r}{\Phi} \\
\node{\A(\star)}            
\arrow{e,t}{W_{\fg}^h}
\node{S(\fg)^{\fg}[[h]]}
\end{diagram}
\hspace{1cm}
\begin{diagram}
\node{\A(\La_{\D})}              
\arrow{e,t}{W_G^h}\arrow{s,r}{\hairnu}
\node{C_{\alg,\D}(G)^G[[h]]} 
\arrow{s,r}{\Phi} \\
\node{\A(\star)}            
\arrow{e,t}{W_{\fg}^h}
\node{S(\fg)^{\fg}[[h]]}
\end{diagram}
$$
Both weight system maps $W_{G}$ and $W_{\fg}$ are defined in terms of
contractions of indices of tensors.  
The commutativity of the left diagram follows from the following
relation between the $\Ad$ and $\ad$ 
$$
\Ad(\exp(x))=\exp(\ad_x)
$$
valid for $x \in \fg$. It remains to deal with the 'wheels factor',
that is with the second diagram. This is dealt by Lemma 
\ref{lem.wheels} below.
\rm{(b)} follows immediately from \rm{(a)} using 
the Definitions of $J_{\fg}$ and $J_G$.
\qed

\begin{lemma}
\lbl{lem.wheels}
The right diagram above commutes.
\end{lemma}

\begin{proof}
We will identify $C_{\alg,\D}(G)^G[[h]]$ with $R_{\D}(T)^W[[h]]$ and
$S(\fg)^{\fg}[[h]]$ with $S(\ft)^W[[h]]$. 
Then, we need to show that the following diagram
$$
\begin{diagram}
\node{\A(\La_{\D})}              
\arrow{e,t}{W_G^h}\arrow{s,r}{\hairnu}
\node{R_{\D}(T)^W[[h]]} 
\arrow{s,r}{\Phi} \\
\node{\A(\star)}            
\arrow{e,t}{W_{\fg}^h}
\node{S(\ft)^W[[h]]}
\end{diagram}
$$
commutes. We will think of $S(\ft)$ as the algebra of polynomials on 
$\ft^\star$. Fix a point $\l \in \ft^\star$ and  a matrix $A 
\in B(\La\to\BZ)$. 
We need to show that for every $\l \in \fg$ we have:
$$
W_{\fg}^h \circ \hairnu (A)(\l)=\Phi \circ W^h_G(A)(\l).
$$
It follows by definition that
$$
W_{\fg}^h \circ \hairnu (A)(\l)=
W_{\fg}^h \circ Z(\text{unknot}) (\l) \cdot 
W_{\fg}^h \circ \hair(A)(\l).
$$
It is well known (see for example \cite[Lemma 2.4]{A}) that
$$
W_{\fg}^h \circ Z(\mathrm{unknot})(\l)=[\l-\rho]=
\frac{j^{1/2}(h \l)}{j^{1/2}(h \rho)}.
$$
Furthermore, we have
$$
\hair(A)=\exp\left(-\frac{1}{2} \sum_{n=1}^\infty a_{2n} \omega_{2n} \right)
$$
where $\omega_{n}$ denotes the wheel with $n$ legs and 
$$
\log (\det(A)(e^h))=\sum_{n=1}^\infty a_{2n} h^{2n} \in \BQ[[h]].
$$
Observe that $W_{\fg}^h: \A(\star)\to S(\fg)^{\fg}[[h]]$ is multiplicative
and the value of $W_{\fg}^h$ on wheels is given by:
$$ 
W_{\fg}^h (w_{2n})(x)=\tr \ad_x^{2n}=2 \sum_{\a \succ 0} \a(x)^{2n}
$$
for $x \in \ft$. Combining this, it follows that
\begin{eqnarray*}
(W_{\fg}^h \circ \hair)(A)(e_{\l}) &=&
W_{\fg}^h \circ 
\exp\left(-\frac{1}{2} \sum_{n=0}^\infty a_{2n} \omega_{2n} \right) \\
&=&
\exp\left(-\frac{1}{2} \sum_{n=0}^\infty a_{2n} W_{\fg}^h(\omega_{2n}) \right)
\\
&=& 
\exp\left(- \sum_{n=0}^\infty a_{2n}  \tr \ad_x^{2n} h^{2n}\right)
\\
&=&
\prod_{\a \succ 0}
\exp\left(- \sum_{n=0}^\infty a_{2n}  (\a,\l)^{2n} h^{2n}\right)
\\
&=&
\prod_{\a \succ 0} 
\exp\left( - \log (\det(A)(e^{h(\a,\l)})) \right) \\
&=&
\prod_{\a \succ 0}
\frac{1}{\det(A)(e^{h(\a,\l)})}.
\end{eqnarray*}
Thus,
$$
W_{\fg}^h \circ \hairnu (A)(\l)=[\l-\rho] 
\prod_{\a \succ 0} \frac{1}{\det(A)(e^{h(\a,\l)})}.
$$
On the other hand, consider the 
element $e^{h\l} \in T$ (where we identify $\ft$ with $\ft^\star$ using
the $W$-invariant inner product) and its Adjoint action $\Ad(e^{h \l})
\in \Aut(\fg)$. It follows that the matrix of $\Ad(e^{h \l})$ is diagonal
with respect to the {\em weight decomposition} of $\fg$:
$$
\fg=\oplus_{a \succ 0} \fg_{-\a} \oplus \ft \oplus
\oplus_{a \succ 0} \fg_{\a}
$$ 
with eigenvalue $e^{h (\a,\l)}$ on each root space $\fg_{\a}$ and
eigenvalue $1$ on the Cartan algebra $\ft$.
Thus, for every Laurrent polynomial $f \in \La$, we have that
$f(\Ad(e^{h \l}))$ is a diagonal matrix with eigenvalue $f(e^{h (\a,\l)})$
on $\fg_{\a}$ and $f(1)$ on $\ft$. Since $A(1)=1$ and $A(t)=A(t^{-1})$,
it follows that 
\begin{eqnarray*}
\Phi \circ W_G^h(A) &=& 
\Phi\left( \frac{1}{\det(\det(A)|_{t\to \Ad(g)}})^{1/2} \right)(\l) \\
& = & [\l -\rho] \frac{1}{\det(\det(A)|_{t \to e^{h(\a,\l)}})^{1/2}} \\
& = & [\l -\rho] \prod_{\a \succ 0} \frac{1}{\det(A)(e^{h(\a,\l)})}.
\end{eqnarray*}
The result follows.
\end{proof}

\ifx\undefined\bysame
	\newcommand{\bysame}{\leavevmode\hbox
to3em{\hrulefill}\,}
\fi

\end{document}